\documentclass[11pt]{article}
\usepackage{graphicx}
\usepackage{amssymb}
\usepackage{amsmath}
\def\n{\noindent}
\def\L{{\bf L}}

\def\ve{\varepsilon}
\def\n{\noindent}
\def\vp{\varphi}

\def\C{{\mathcal C}}

\def\S{{\mathcal S}}

\def\diams{\diamondsuit}
\def\ds{\displaystyle}
\def\sqr#1#2{\vbox{\hrule height .#2pt
\hbox{\vrule width .#2pt height #1pt \kern #1pt
\vrule width .#2pt}\hrule height .#2pt }}
\def\square{\sqr74}
\def\endproof{\hphantom{MM}\hfill\llap{$\square$}\goodbreak}

\def\bega{\begin{array}}
\def\enda{\end{array}}
\def\begi{\begin{itemize}}
\def\endi{\end{itemize}}

\def\P{{\mathcal P}}

\def\N{{\cal N}}

\def\R{I\!\!R}

\def\v{\vskip 1em}
\def\vs{\vskip 2em}

\def\be{\begin{equation}}
\def\beq{\begin{equation}}
\def\bel{\begin{equation}\label}
\def\eeq{\end{equation}}

\begin{document}
\title{\bf Uniqueness of Conservative Solutions to the Camassa-Holm Equation via Characteristics}
\vs
\author{Alberto Bressan$^{(*)}$, Geng Chen$^{(**)}$, and Qingtian Zhang$^{(*)}$\\    \\
{\small (*) Department of Mathematics, Penn State University,
University Park, Pa.~16802, U.S.A.}\\
{\small (**) School of Mathematics,
Georgia Institute of Technology,
Atlanta, Ga.~30332, U.S.A.}
\\
\,
\\
{\small e-mails:~ bressan@math.psu.edu,~gchen73@math.gatech.edu, ~zhang\_q@math.psu.edu}}

\maketitle

\begin{abstract} The paper provides a direct proof the uniqueness
of  solutions to the Camassa-Holm equation, based on characteristics.
Given a conservative solution $u=u(t,x)$,
an equation is introduced which singles out a unique characteristic
curve through each initial point.
By studying the evolution of the quantities $u$ and $v= 2\arctan u_x$
along each characteristic, it is proved that the Cauchy problem with general
initial data
$u_0\in H^1(\R)$ has a unique solution, globally in time.
\end{abstract}

\v
\section{Introduction}
\setcounter{equation}{0}
\v
The Cauchy problem for the Camassa-Holm equation \cite{CH} can be written in the form
\bel{1.1}u_t+(u^2/2)_x+P_x~=~0\,,\eeq
\bel{id}u(0,x)~=~u_0(x).\eeq
The nonlocal source term $P$ is here defined as a convolution:
\bel{Pdef}P~\doteq~{1\over 2} e^{-|x|} * \left(u^2+{u_x^2\over 2}
\right)\,.\eeq
For any initial data $u_0\in H^1(\R)$, various papers have studied
the global existence of solutions \cite{BC2, BC3, HR, XZ1, XZ2}.  In particular,
in \cite{BC2} a new set of variables was introduced, transforming the
equation \ref{1.1} into a semilinear system.   This yields a group of
conservative (i.e., energy-preserving) solutions $u(t) = S_t u_0$, depending
continuously on the initial data w.r.t.~suitable norms.
 In \cite{BC3} a similar approach was used to
construct a semigroup of dissipative solutions.
Solutions obtained by this particular
transformation of variables are clearly unique.  However, in principle, other constructive procedures may
yield different solutions to the same Cauchy problem.

Uniqueness  is a delicate issue because in general the flow  map $(t, u_0)
\mapsto S_t u_0$ constructed in \cite{BC2}
is not continuous as a map from $[0,T]\times H^1$ into $H^1$,
neither as a function of time nor of the initial data.  Rather, it is continuous from
$[0,T]\times H^1$ into spaces with weaker norms such as $\L^2$ or $\C^0$.

The papers \cite{BF, GHR1, GHR2}
have introduced new distances $d^\diams(\cdot, \cdot)$ on $H^1$ which render the
flow map uniformly Lipschitz continuous on bounded subsets of $H^1$. These distances
satisfy
\bel{Lipd}
{d\over dt} \Big[d^\diams(S_t u_0, S_t \tilde u_0)
\Big]~\leq~C\,d^\diams(S_t u_0, S_t \tilde u_0),
\eeq
for every pair of initial data $u_0,\tilde u_0\in H^1$
and some constant $C$ depending only on the $H^1$ norm of $u_0, \tilde u_0$.
We can now consider any map $t\mapsto w(t) $ from an interval $[0,T]$ into $H^1$ which is Lipschitz
continuous w.r.t.~the distance
$d^\diams(\cdot, \cdot)$ and satisfies
\bel{tan}\liminf_{h\to 0+} {1\over h}\,d^\diams
\Big(w(t+h)\, , ~ S_h w(t)
\Big)~=~0\qquad
 \hbox{for a.e.}~t\geq 0.\eeq
An elementary error estimate (see for example \cite{Bbook}, p.~26) then
yields $w(t) = S_t w(0)$ for every $t\ge 0$.  This approach
provides some result on
the uniqueness of solutions.  In practice, however, checking that any conservative
solution of (\ref{1.1}) satisfies the tangency condition
(\ref{tan}) is not an easy task.

\begin{figure}[htbp]
\centering
  \includegraphics[scale=0.550]{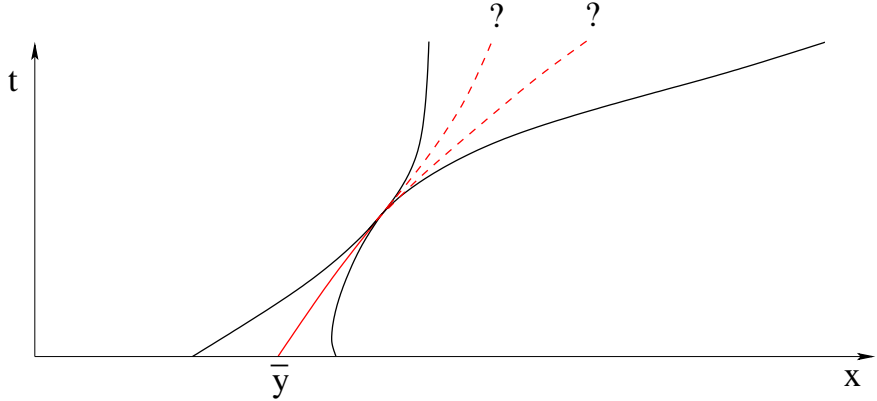}
    \caption{\small When $u$ is only H\"older continuous, the equation (\ref{car})
    does not determine a unique characteristic
    starting at the point $\bar y$.}
\label{f:wa11}
\end{figure}

Motivated by a recent paper by Dafermos \cite{Daf},
our present goal is to develop a direct approach to the uniqueness of conservative solutions, based on characteristics.
The underlying idea is simply to write a set of ODEs satisfied by $u$ and $u_x$ along
the characteristic starting at a given point $\bar y$.
  If this set of equations has a unique solution
for a.e.~$\bar y$, then the entire solution $u(\cdot)$ will be uniquely determined.
At this naive level, the approach runs into a fundamental
difficulty.  Namely,  since the solution $u$ is only H\"older continuous,
the Cauchy problem
\bel{car}{d\over dt} x(t)~=~u(t, x(t)),\qquad\qquad x(0) ~=~\bar y\eeq
may well have several solutions   (Fig.~\ref{f:wa11}).
To overcome this stumbling block, our analysis relies on two key ideas.
\begi
\item{} Since we assume that $u$ is a conservative solution, the quantity
 $w=u_x^2$ provides a measure-valued solution to  the balance law
\bel{blaw}w_t+(uw)_x
~=~2(u^2-P) u_x\,. \eeq
Because of (\ref{blaw}),  the
characteristic curve $t\mapsto x(t)$   satisfies the additional
equation
\bel{ceq}
{d\over dt} \int_{-\infty}^{x(t)}u_x^2(t,x)\, dx~=~
\int_{-\infty}^{x(t)}[2u^2 u_x - 2Pu_x  ](t,x)\, dx\,,\qquad\qquad x(0)=\bar y.
\eeq
By itself, (\ref{ceq}) still does not single out a unique characteristic
(think for example of the case where $u(t,x)\equiv 0$).   However, combining
 two equations  (\ref{car}) and (\ref{ceq}) we eventually
obtain an integral equation with unique solutions.

\item{}  Instead of the variables $(t,x)$, it is convenient to work with
an adapted set of variables $(t, \beta)$, where
$\beta$ is implicitly defined as
\bel{beta} x(t,\beta) + \int_{-\infty}^{x(t,\beta)} u_x^2(t,\xi)\, d\xi~=~\beta\,.\eeq
In terms of these variables, the solution $u=u(t,\beta)$ becomes globally Lipschitz continuous.
Indeed, $|u_\beta|\leq 1$ while $|u_t|\leq C$ for some constant $C$ depending only on
the $H^1$ norm of the solution.
\endi
Our analysis eventually shows that, for any conservative solution to the Camassa-Holm equation (\ref{1.1}), the characteristic curves $t\mapsto x(t,\bar y)$
as well as the values of $u$ and $v\doteq 2\arctan u_x$ along these curves
can be recovered
by a set of integro-differential equations
having unique solutions.   In turn, this provides a direct proof
of the uniqueness of conservative solutions to (\ref{1.1}),
for general initial data $u_0\in H^1(\R)$.

We emphasize the  major difference between the present approach  and
previous ones.
In \cite{BC2} one starts by constructing a solution to an auxiliary semilinear system.
After a suitable variable transformation, this yields a conservative solution
to the Camassa-Holm equation (\ref{1.1}).
Here we follow an inverse route.  Given a conservative
solution
$u=u(t,x)$ to (\ref{1.1}), we introduce a set of auxiliary variables
tailored to this particular
solution.   We then prove that these variables satisfy a particular semilinear system
having unique solutions. In turn, this yields the uniqueness of the
conservative solution $u$ in the original variables.

The remainder of the paper is organized as follows.  In  Section~2 we review
basic definitions and state our main uniqueness result, Theorem~2.
Section~3 establishes the key technical tool (Lemma~2), determining a unique
characteristic curve through each initial point. In Section~4 we study how the
gradient $u_x$ of a conservative solution varies along a characteristic, and
conclude the proof of the main theorem.

\section{Basic definitions and results}
\setcounter{equation}{0}
To make sense of the source term $P$, at each time $t$
we require that the function $u(t,\cdot)$
lies in the space $H^1(\R)$ of absolutely
continuous
functions $u\in\L^2(\R)$ with derivative $u_x\in\L^2(\R)$, endowed
with the norm
$$\big\| u \big\|_{H^1}~\doteq~
\Bigl( \int_{\R} \big[u^2(x)+u^2_x(x)\big]\,dx\Bigr)^{1 / 2}.$$
For $u \in H^1(\R)$, Young's inequality ensures that
$$P~=~(1-\partial_x^2)^{-1}\Big(u^2+{u_x^2 \over 2}\Big)~ \in~ H^1(\R).$$
For future use we record the following inequalities, valid for any function $u\in H^1(\R)$:
\bel{uli}
\|u\|_{\L^\infty}~\leq~\|u\|_{H^1}\,,\eeq
\bel{1.8}\| P\|_{\L^\infty}\,,~~\|
P_x\|_{\L^\infty}~~\leq
~~\Big\|{1\over 2} e^{-|x|}\Big\|_{\L^\infty} \cdot \left\|
u^2+{u_x^2\over 2}
\right\|_{\L^1}~~\leq~~{1\over 2} \,\|u\|_{H^1}\,,\eeq
\bel{1.9}\| P\|_{\L^2}\,,~~\| P_x\|_{\L^2}~~\leq
~~\Big\|{1\over 2} e^{-|x|}\Big\|_{\L^2} \cdot \left\| u^2+{u_x^2\over
2}
\right\|_{\L^1}~~\leq~~{1\over \sqrt 2} \,\|u\|_{H^1}\,.\eeq

\v
\n{\bf Definition 1.}
{\it By a solution of the Cauchy problem (\ref{1.1})-(\ref{id}) on $[0,T]$
we mean
a H\"older continuous function $u=u(t,x)$ defined on
$[0,T]\times\R$
with the
following properties. At each fixed $t$ we have $u(t,\cdot)
\in H^1(\R)$. Moreover,
the map $t\mapsto u(t,\cdot)$ is Lipschitz continuous
from $[0,T]$ into $\L^2(\R)$ and satisfies the initial condition
(\ref{id})
together with
\bel{1.3}{d\over dt} \,u ~= ~- uu_x-P_x\eeq
for a.e.~$t$.  Here (\ref{1.3}) is understood as an equality
between functions in $\L^2(\R)$.}
\v
As shown in \cite{BC2, BC3}, as soon as the gradient of a solution blows up,
uniqueness is lost, in general.  To single out a unique solution, some additional
conditions are needed.

For smooth solutions, differentiating (\ref{1.1}) w.r.t.~$x$ one obtains
\bel{1.4}u_{xt}+(uu_x)_x~=~\left( u^2+{u_x^2\over 2}\right)-P\,.\eeq
Multiplying (\ref{1.1}) by $u$ and (\ref{1.4}) by $u_x$, we
obtain the two conservation laws with source term
\bel{1.5}\left(u^2\over 2\right)_t+\left({u^3\over 3}+u\,P\right)_x~=~u_x
\,P\,,
\eeq
\bel{1.6}\left(u^2_x\over 2\right)_t+\left({uu_x^2\over 2}-{u^3\over
3}\right)_x
~=~- u_x \,P. \eeq
Summing (\ref{1.5}) and (\ref{1.6}), and integrating w.r.t.~$x$, we see that
for smooth solutions
the total energy
\bel{1.7}E(t)~\doteq~ \int_{\R} \Big( u^2(t,x)+u_x^2(t,x)\Big)\,dx\eeq
is constant in time.
\v

\n{\bf Definition 2.}
{\it A solution $u=u(t,x)$ is {\em conservative} if $w=u_x^2$ provides
a distributional solution to the balance law (\ref{blaw}), namely
\bel{ec}
\int_0^\infty\int \Big[u_x^2\vp_t+uu_x^2\vp_x + 2\,(u^2-P)\,u_x \vp
\Big] \, dx dt+\int u_{0,x}^2(x)
 \vp(0,x)\, dx~=~0\eeq
for every test function $\vp\in \C^1_c(\R^2)$.}
\v
 The main result proved in
\cite{BC2}, on the global existence of conservative solutions can be stated as follows.
\v
\n{\bf Theorem 1.} {\it For any initial data $u_0\in H^1(\R)$ the Camassa-Holm equation has
a global conservative solution $u=u(t,x)$.
More precisely, there exists a family of Radon measures
$\{\mu_{(t)}\,,~t\in\R\}$, depending continuously on time w.r.t.~the topology of weak convergence of measures, such that the following
properties hold.
\v
\begi
\item[(i)]
The function $u$ provides a solution to the Cauchy problem
(\ref{1.1})-(\ref{id}) in the sense of Definition~1.
\v
\item[(ii)] There exists a null set $\N\subset\R$ with $meas(\N)=0$
such that for every $t\notin \N$ the measure $\mu_{(t)}$ is absolutely
continuous and has density $u^2_x(t,\cdot)$
w.r.t.~Lebesgue measure.
\v
\item[(iii)]
The
family $\{\mu_{(t)}\,;~t\in\R\}$
provides a measure-valued solution $w$ to the linear transport equation
with source
\bel{5.20}w_t+(uw)_x~=~2\,(u^2-P)\,u_x\,.\eeq
\endi}

At a time  $t\in\N$  the measure
$\mu_{(t)}$ has a nontrivial singular part.  For a conservative solution
$u$ which is not smooth, in general we only know that the
energy $E$ in   (\ref{1.7}) coincides a.e.~with a constant.   Namely,
$$E(t) ~=~E(0)\quad\hbox{for}~~t\notin \N,\qquad\qquad E(t) ~<~E(0)
\quad\hbox{for}~~t\in \N.$$

The main purpose of this paper is to prove the uniqueness of the above solution.
\v
{\bf Theorem 2.} {\it For any initial data $u_0\in H^1(\R)$,
the Cauchy problem (\ref{1.1})-(\ref{id}) has a unique conservative solution.
 }
 \v

\section{Preliminary lemmas}
\setcounter{equation}{0}

Let $u=u(t,x)$ be a solution of the Cauchy problem (\ref{1.1})-(\ref{id})
which satisfies the additional balance law (\ref{ec}).
As mentioned in the Introduction, it is  convenient to work with
the adapted coordinates $(t,\beta)$, related to the original coordinates
$(t,x)$ by the integral relation (\ref{beta}).
At  times $t$ where the measure $\mu_{(t)}$
is not absolutely continuous w.r.t.~Lebesgue measure,  we can  define
$x(t,\beta)$
to be the unique point $x$ such that
\bel{xgen}
x(t,\beta)+\mu_{(t)}( \,]-\infty, x[\,)~\leq~ \beta~\leq~x(t,\beta)
+\mu_{(t)}(\, ]-\infty, x]).
\eeq
Notice that (\ref{xgen}) and (\ref{beta}) coincide at every time where $
\mu_{(t)}$ is absolutely continuous with density $u_x^2$ w.r.t.~Lebesgue measure.
The next lemma, together with Lemma~3,
establishes the Lipschitz continuity of $x$ and $u$ as functions
of the variables $t,\beta$.
\v
{\bf Lemma 1.} {\it Let $u=u(t,x)$ be a conservative solution of (\ref{1.1}).
Then, for every $t\ge 0$,  the  maps
$\beta\mapsto x(t,\beta)$ and $\beta\mapsto u(t,\beta)\doteq u(t, x(t, \beta))$
implicitly defined by (\ref{xgen})
are Lipschitz continuous with constant $1$.  The map $t\mapsto x(t,\beta)$
is also Lipschitz continuous with a constant depending only on $\|u_0\|_{H^1}$.
}

{\bf Proof.}
{\bf 1.} Fix any time $t\geq 0$.  The
the map
$$x~\mapsto ~\beta(t,x)~\doteq ~x+\int_{-\infty}^x u_x^2(t,y)\, dy$$
is right continuous and strictly increasing.  Hence it has a well defined, continuous,
nondecreasing
inverse $\beta\mapsto x(t,\beta)$.
If  $\beta_1<\beta_2$, then
\bel{b12}
x(t,\beta_2)- \-x(t,\beta_1)+\mu_{(t)}\Big( ]x(t,\beta_1)\,,~x(t,\beta_2)[\Big ) ~\le~
\beta_2-\beta_1\,.
\eeq
This implies
$$x(t,\beta_2)-x(t,\beta_1)~\leq~ \beta_2-\beta_1\,,$$
showing that the map $\beta\mapsto x(t,\beta)$ is a contraction.
\v
{\bf 2.} To prove the
 Lipschitz continuity of the map $\beta\mapsto u(t,\beta)$, assume $\beta_1<\beta_2$.
By (\ref{b12}) it follows
\bel{ulip}\bega{c}\ds
\Big|u(t,x(t,\beta_2))-u(t,x(t,\beta_1))\Big|~\leq~
\int_{x(t,\beta_1)}^{x(t,\beta_2)}|u_x| dx~\leq~
 \int_{x(t,\beta_1)}^{x(t,\beta_2)}\frac12(1+ u_x^2)  dx\cr\cr
 \ds \leq~  \frac12\left[x(t,\beta_2)- \-x(t,\beta_1)+\mu_{(t)}\Big( ]x(t,\beta_1)\,,~x(t,\beta_2)[\Big )\right] ~\le~
 \frac12(\beta_2-\beta_1)\,.\enda
\eeq
\v
{\bf 3.} Next, we prove the Lipschitz continuity of the map
$t\mapsto x(t,\beta)$.    Assume $x(\tau,\beta) =y$.
We recall that the family of measures $\mu_{(t)}$ satisfies the balance law (\ref{5.20}),
where for each $t$ the drift $u$ and the source term  $2(u^2-P)u_x$ satisfy
\bel{uib}
\|u\|_{\L^\infty(\R)}~\leq~C_\infty~\doteq~\|u\|_{H^1(\R)}\,,\eeq
\bel{stb}\|2(u^2-P)u_x\|_{\L^1(\R)}~\leq~2\Big( \|u\|_{\L^\infty}\|u\|_{\L^2}+\|P\|_{\L^2}
\Big) \|u_x\|_{\L^2}~\leq~C_S\,,\eeq
for some constant $C_S$ depending only on the $H^1$ norm of the solution.
Therefore, for $t>\tau$ we have
$$\bega{rl}\mu_{(t)}\Big(]-\infty\,,~ y-C_\infty(t-\tau)[\Big)&\ds\leq~
\mu_{(\tau)}\Big(]-\infty\,,~ y[\Big)+\int_\tau^t \|2(u^2-P)u_x\|_{\L^1(\R)}\, dt\cr\cr
&\leq~\mu_{(\tau)}\Big(]-\infty\,,~ y[\Big)+C_S(t-\tau)\,. \enda$$
Defining $y^-(t)\doteq y-(C_\infty+C_S)(t-\tau)$, we obtain
$$\bega{l}y^-(t)+\mu_{(t)}\Big(]-\infty\,,~ y^-(t)[\Big)~\leq~
y-(C_\infty+C_S)(t-\tau)+\mu_{(\tau)}\Big(]-\infty\,,~ y[\Big)+C_S(t-\tau)\cr\cr
\qquad \le~y+\mu_{(\tau)}\Big(]-\infty\,,~ y[\Big)~\leq~\beta\,.\enda$$
This implies $x(t,\beta)~\geq~ y^-(t)$ for all $t>\tau$.    An entirely similar
argument yields $x(t,\beta)~\leq~ y^+(t)\doteq y+(C_\infty+C_S)(t-\tau)$,
proving the uniform Lipschitz continuity of the map $t\mapsto x(t,\beta)$.
 \endproof
\v
The next result, which provides the foundation to all our analysis,
shows that characteristics can be uniquely determined by an integral equation
combining (\ref{car}) with (\ref{ceq}).

{\bf Lemma 2.} {\it Let $u=u(t,x)$ be a conservative solution of the Camassa-Holm
equation. Then, for any $\bar y\in\R$ there exists a unique
Lipschitz continuous map  $t\mapsto x(t)$ which
satisfies both  (\ref{car}) and (\ref{ceq}).
In addition, for any $0\le \tau\leq t$ one has
\bel{ucar}
u(t, x(t)) - u(\tau, x(\tau)) ~=~-\int_\tau^t P_x(s, x(s))\, ds\,.\eeq
}

{\bf Proof.} {\bf 1.} Using the adapted  coordinates $(t, \beta)$ as in (\ref{beta}), we
write
the characteristic
starting at $\bar y$ in the form $t\mapsto x(t)=x(t, \beta(t))$, where
$\beta(\cdot)$ is a map to be determined.
By summing the two equations (\ref{car}) and (\ref{ceq}) and integrating w.r.t.~time
we obtain
\bel{ieq}x(t) + \int_{-\infty}^{x(t)}u_x^2(t,x)\, dx~=~\bar y + \int_{-\infty}^{\bar y}
u_{0,x}^2(x)\, dx + \int_0^t \left( u(t, x(t)) + \int_{-\infty}^{x(t)}[2u^2 u_x - 2Pu_x  ](t,x)\, dx\right)dt\,.
\eeq
Introducing the
function
\bel{Gdef}
G(t, \beta)~ \doteq~
\int_{-\infty}^{x(t,\beta)}[ u_x + 2u^2u_x- 2u_x P]
\, dx\eeq
and the constant
\bel{bbeta}\bar\beta~=~\bar y + \int_{-\infty}^{\bar y}
u_{0,x}^2(x)\, dx\,,\eeq
we can rewrite the equation (\ref{ieq}) in the form
\bel{bt}
\beta(t)~=~\bar\beta+\int_0^t G(s,\beta(s))\, ds\,.\eeq
\v
{\bf 2.} For each fixed $t\geq 0$,
since the maps $x\mapsto u(t,x)$ and $x\mapsto   P(t,x)$
are both in $ H^1(\R)$,  the function $\beta\mapsto G(t,\beta)$
defined at (\ref{Gdef}) is
uniformly bounded and absolutely continuous.
Moreover,
\bel{Glip}G_\beta~=~[u_x + 2u^2 u_x  - 2u_x  P ]\, x_\beta~=~
{u_x + 2u^2 u_x  - 2u_x  P\over 1+u_x^2}~ \in ~[-C,\,C]\eeq
for some constant $C$ depending only on the $H^1$ norm of $u$.
Hence the function $G$ in (\ref{Gdef}) is uniformly Lipschitz continuous w.r.t.~$\beta$.
\v
{\bf 3.} Thanks to the Lipschitz continuity of the function $G$, the existence of a unique
solution to the
integral equation (\ref{bt})
can be proved by a standard fixed point argument.  Namely, consider the
Banach space of all continuous functions
$\beta:\R_+\mapsto \R$ with weighted norm
$$\|\beta\|_*~\doteq~\sup_{t\geq 0} ~e^{-2Ct} |\beta(t)|\,.$$
On this space, we claim that the Picard map
$$(\P \beta)(t)~\doteq~\bar \beta+\int_0^t G(\tau, \beta(\tau))\, d\tau$$
is a strict contraction. Indeed, assume $\|\beta-\tilde \beta\|_*=\delta>0$.
This implies
$$|\beta(\tau)-\tilde\beta(\tau)|~\leq ~\delta e^{2C\tau}\qquad\qquad
\hbox{for all}~\tau\geq 0.$$
Hence
$$\bega{rl}\Big|(\P\beta)(t) - (\P\tilde\beta)(t)\Big|&\ds
=~\left|\int_0^t \Big(G(\tau, \beta(\tau))-
G(\tau, \tilde \beta(\tau))\Big)\, d\tau\right|~\leq~
\int_0^t C\Big|\beta(\tau)-
\tilde \beta(\tau)\Big|\, d\tau\cr\cr
&\ds\leq~\int_0^t C\delta e^{2C\tau}\, d\tau~\leq~{\delta\over 2} \, e^{2Ct}.
\enda
$$
We thus conclude that $\|\P\beta-\P\tilde \beta\|_*\leq \delta/2$.

By the contraction mapping principle, the integral equation (\ref{bt}) has a unique
solution.
 \v
\begin{figure}[htbp]
\centering
  \includegraphics[scale=0.550]{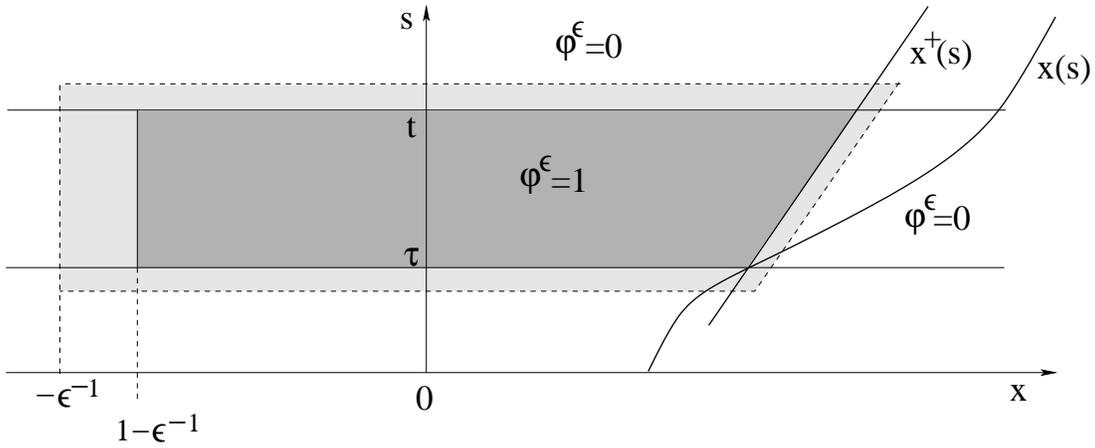}
    \caption{\small The Lipschitz continuous test function $\vp^{\epsilon}$
introduced at (\ref{pen1}). }
\label{f:wa15}
\end{figure}

{\bf 4.} By the previous construction, the map
$t\mapsto x(t)\doteq x(t, \beta(t))$ provides the unique solution to (\ref{ieq}).
Being the composition of two Lipschitz  functions, the map $t\mapsto x(t,\beta(t))$
is Lipschitz continuous.
To prove that it
 satisfies the ODE for characteristics (\ref{car}), it suffices to show that
 (\ref{car}) holds at each time $\tau>0$ such that
 \begi
 \item[(i)] $x(\cdot)$ is differentiable at $t=\tau$, and
\item[(ii)]   the measure $\mu_{(\tau)}$ is absolutely continuous.
\endi

Assume, on the contrary, that $\dot x(\tau)\not= u(\tau, x(\tau))$.
To fix the ideas, let
 \bel{ass}\dot x(\tau)~=~ u(\tau, x(\tau))+2\ve_0\eeq
 for some $\ve_0>0$.   The case $\ve_0 <0$ is entirely similar.
To derive a contradiction we
observe that, for all $t\in ]\tau, \tau+\delta]$,
with $\delta>0$ small enough one has
\bel{xpmt}x^+(t)~\doteq~
x(\tau) +(t-\tau) [u(\tau, x(\tau))+\ve_0 ]~<~x(t)\,.\eeq
We also observe that, since $u,P$ are continuous while $u_x\in\L^2$, by an
approximation argument the identity in (\ref{ec}) remains
valid for any test function $\vp\in H^1$ with compact support.
In particular, this is true if $\vp$ is
Lipschitz continuous with compact support.

For any $\epsilon>0$ small, we can thus consider the functions
$$\rho^{\epsilon}(s,y)~\doteq~\left\{\bega{cl} 0 \qquad &\hbox{if}\quad y
\leq -\epsilon^{-1},\cr
(y+\epsilon^{-1}) \qquad &\hbox{if}\quad -\epsilon^{-1}\leq y\leq 1-\epsilon^{-1},\cr
1\qquad &\hbox{if}\quad 1-\epsilon^{-1} \leq y\leq x^+(s),\cr
1-\epsilon^{-1}(y-x(s))\qquad &\hbox{if}\quad  x^+(s)\leq y\leq x(s)^++\epsilon,\cr
0\qquad &\hbox{if}\quad y\geq x^+(s)+\epsilon,\enda\right.$$
\bel{timtest}\chi^\epsilon(s)~\doteq~\left\{\bega{cl} 0\qquad 
&\hbox{if}\quad s\leq \tau-\epsilon,\cr
\epsilon^{-1}(s-\tau+\epsilon)\qquad &\hbox{if}\quad \tau-\epsilon\leq s\leq \tau,\cr
1\qquad &\hbox{if}\quad \tau\leq s\leq  t,\cr
1-\epsilon^{-1}(s-t) \qquad &\hbox{if}\quad t\leq s<t+\epsilon,\cr
0 \qquad &\hbox{if}\quad s\geq t+\epsilon.\enda\right.
\eeq

Define \bel{pen1}\varphi^{\epsilon}(s,y)~\doteq~\min\{ \varrho^{\epsilon}(s,y),
\,\chi^\epsilon(s)\}.\eeq

Using $\vp^{\epsilon}$ as test function  in (\ref{ec}) we obtain
\bel{vpe}
\iint \Big[u_x^2\vp^{\epsilon}_t+uu_x^2\vp^{\epsilon}_x 
+ 2\,(u^2-P)\,u_x \vp^{\epsilon}
\Big] \, dx dt~=~0.
\eeq

We now observe that, if $t$ is sufficiently close to $\tau$, then
\bel{bdt}\lim\limits_{\epsilon\to 0}\int_\tau^t
\int_{x^+(s)-\epsilon}^{x^+(s)+\epsilon}
u_x^2(\vp^{\epsilon}_t+u\vp^{\epsilon}_x)\, dx ds
~\geq~0\eeq
Indeed, for $s\in [\tau+\epsilon, \, t-\epsilon]$ one has
$$0 ~=~ \vp^{\epsilon}_t + [u(\tau, x(\tau))+\ve_0 ] \vp^{\epsilon}_x ~
\leq~\vp^{\epsilon}_t + u(s,x) \vp^{\epsilon}_x\,,$$
because  $u(s,x)<u(\tau, x(\tau))+\ve_0 $ and $\vp^{\epsilon}_x\leq 0$.

Since the family of measures $\mu_{(t)}$ depends continuously on $t$
in the topology of weak convergence, taking the limit of (\ref{vpe}) as
 $\epsilon\to 0$, for $\tau,t\notin\N$ we obtain
\bel{55}\bega{rl} 0&=~
\ds\int_{-\infty}^{x(\tau)}u_x^2(\tau,x)\, dx-\int_{-\infty}^{x^+(t)}u_x^2(t,x)\, dx
+ \int_\tau^t\int_{-\infty}^{x^+(s)}[2u^2 u_x - 2u_x P ]\,dxds\cr\cr
\ds &\qquad\qquad \ds +\lim\limits_{\epsilon\to 0}\int_\tau^t
\int_{x^+(s)-\epsilon}^{x^+(s)+\epsilon}  u_x^2(\vp^\epsilon_t+u\vp^\epsilon_x)\, dx ds
\cr\cr
& \ds
\geq~
\ds\int_{-\infty}^{x(\tau)}u_x^2(\tau,x)\, dx-\int_{-\infty}^{x^+(t)}u_x^2(t,x)\, dx
+ \int_\tau^t\int_{-\infty}^{x^+(s)}[2u^2 u_x - 2u_x P ]\,dxds
\enda\eeq
In turn, (\ref{55})  implies
$$\bega{l}\ds \int_{-\infty}^{x^+(t)}u_x^2(t,x)\, dx~\geq~
\int_{-\infty}^{x(\tau)}u_x^2(\tau,x)\, dx
+ \int_\tau^t\int_{-\infty}^{x^+(s)}[2u^2 u_x - 2u_x P ]\,dxds\cr\cr
\ds
=~\ds\int_{-\infty}^{x^+(t)}u_x^2(t,x)\, dx
+ \int_\tau^t\int_{-\infty}^{x(s)}[2u^2 u_x - 2u_x P ]\,dxds+ o_1(t-\tau).
\enda
$$
Notice that the last term is a higher order infinitesimal,
 satisfying
${o_1(t-\tau)\over t-\tau}\to 0$ as $t\to \tau$.
Indeed
$$\bega{l}\ds|o_1(t-\tau)|~=~\left|\int_\tau^t \int_{x^+(s)}^{x(s)}[2u^2y_x - 2 Pu_x]\,dx ds\right|~
\leq~\int_\tau^t \int_{x^+(s)}^{x(s)}|2u^2 - 2 P|\, |u_x|\,dxds\cr\cr
\qquad\ds\leq~2\|u^2-P\|_{\L^\infty} \cdot \int_\tau^t \int_{x^+(s)}^{x(s)}\, |u_x|\,dxds
\cr\cr\qquad\ds\leq~2\|u^2-P\|_{\L^\infty} \int_\tau^t (x(s)-x^+(s))^{1/2}\, \|u_x(s,\cdot)\|_{\L^2}\,ds~ \leq~C\cdot (t-\tau)^{3/2}.\enda$$

On the other hand, by (\ref{Gdef}) and (\ref{bt}) a linear approximation yields
\bel{3.18}\beta(t) ~=~\beta(\tau) +(t-\tau) \left[ u(\tau, x(\tau)) + {2\over 3} u^3(\tau, x(\tau))
 -\int_{-\infty}^{x(\tau)}2u_xP\, dx\right]+ o_2(t-\tau)\,,\eeq
with ${o_2(t-\tau)\over t-\tau}\to 0$ as $t\to \tau$.

For every $t>\tau$ with $t\notin \N$, $t$ sufficiently close to $\tau$, we now have
\bel{3.19}
\bega{lll}
\ds\beta(t)&\ds=&\ds x(t)+\int_{-\infty}^{x(t)}u_x^2(t,x)dx\cr\cr
&\ds >&\ds x(\tau)+(t-\tau)[u(\tau,x(\tau))+\ve_0 ]+\int_{-\infty}^{x^+(t)}u_x^2(t,y)dy\cr\cr
&\ds \geq&\ds x(\tau)+(t-\tau)[u(\tau,x(\tau))+\ve_0 ]+
\int_{-\infty}^{x(\tau)}u_x^2(\tau,x)dx\cr\cr
&&\qquad \ds +\int_\tau^t\int_{-\infty}^{x(s)} [2u^2 u_x - 2u_x P ]\, dxds+ o_1(t-\tau).
\enda
\eeq
Combining (\ref{3.18}) and (\ref{3.19}), we find
\beq\bega{l}
\ds\beta(\tau) +(t-\tau) \left[ u(\tau, x(\tau)) +\int_{-\infty}^{x(\tau)}[2u^2 u_x - 2u_x P ]\, dx\right]+o_2(t-\tau)\cr\cr
\qquad \ds\geq ~\left[x(\tau)+\int_{-\infty}^{x(\tau)}u_x^2(\tau,x)\,dx\right]
+(t-\tau)[u(\tau,x(\tau))+\ve_0 ]\cr\cr
\qquad\qquad \ds+\int_\tau^t\int_{-\infty}^{x(s)} [2u^2 u_x - 2u_x P ]\, dxds
+o_1(t-\tau).
\enda
\eeq
Subtracting common terms, dividing both sides by $t-\tau$ and
letting $t\to \tau$, we achieve a contradiction.  Therefore, (\ref{car}) must hold.
\v
{\bf 5.}  We now prove (\ref{ucar}).
By (\ref{1.3}), for every test function $\phi\in \C^\infty_c(\R^2)$ one has
\bel{CHW1}
\int_0^\infty \int\left[u\phi_t +{u^2\over 2}\phi_x + P_x\phi\right]\, dxdt
+ \int u_0(x)\phi(0,x)\, dx~=~0\,.\eeq
Given any $\vp\in\C^\infty_c$, let $\phi=\vp_x$.
Since the map $x\mapsto u(t,x)$ is absolutely continuous, we can
integrate by parts w.r.t.~$x$ and obtain
\bel{CHW2}
\int_0^\infty \int\left[u_x\vp_t +uu_x\vp_x + P_x\vp_x\right]\, dxdt
+ \int u_{0,x}(x)\vp(0,x)\, dx~=~0\,.\eeq

By an
approximation argument, the identity (\ref{CHW2}) remains
valid for any test function $\vp$ which is
Lipschitz continuous with compact support.
For any $\epsilon>0$ sufficiently small, we thus consider the functions 
$$\varrho^{\epsilon}(s,y)~\doteq~\left\{\bega{cl} 0 \qquad &\hbox{if}\quad y\leq -
\epsilon^{-1},\cr
y+\epsilon^{-1} \qquad &\hbox{if}\quad -\epsilon^{-1}\leq y\leq 1-\epsilon^{-1},\cr
1\qquad &\hbox{if}\quad 1-\epsilon^{-1} \leq y\leq x(s),\cr
1-\epsilon^{-1}(y-x(s))\qquad &\hbox{if}\quad  x(s)\leq y\leq x(s)+\epsilon,\cr
0\qquad &\hbox{if}\quad y\geq x(s)+\epsilon,\enda\right.$$

\bel{pen}\psi^{\epsilon}(s,y)~\doteq~\min\{ \varrho^{\epsilon}(s,y),\,\chi^\epsilon(s)\},\eeq
where $\chi^\epsilon(s)$ as in (\ref{timtest}).
We now use the test function
$\vp=\psi^{\epsilon}$  in (\ref{CHW2}) and let
$\epsilon\to 0$. Observing that the function $P_x$
is continuous, we obtain
\bel{Inteq}\bega{rl}
\ds \int_{-\infty}^{x(t)} u_x(t,x)\, dx&=~\ds\int_{-\infty}^{x(\tau)}u_x(\tau,x)\, dx
- \int_\tau^t P_x(s, x(s))\,ds\cr\cr
&\ds+\lim\limits_{\epsilon\to 0}\int_{\tau-\epsilon}^{t+\epsilon}\int_{x(s)}
^{x(s)+\epsilon} u_x(\psi^{\epsilon}_t+u\psi^{\epsilon}_x) dx ds\,.\cr\cr
\enda\eeq
To complete the proof it suffices to show that the last term on the right hand side of
(\ref{Inteq}) vanishes. Since $u_x\in \L^2$, Cauchy's inequality yields
\bel{lexx}
\bega{ll}
\ds\left|\int_{\tau}^{t}\int_{x(s)}
^{x(s)+\epsilon} u_x(\psi^{\epsilon}_t+u\psi^{\epsilon}_x) dx ds\right|
&\ds\leq\int_{\tau}^{t}\left(\int_{x(s)}^{x(s)+\epsilon}
 |u_x|^2dx\right)^{1/2}\left(\int_{x(s)}^{x(s)+\epsilon}
(\psi^\epsilon_t+u\psi^\epsilon_x)^2dx\right)^{1/2}ds\,.\cr\cr
\enda
\eeq
For each $\epsilon>0$
consider the function
\bel{ete}\eta_\epsilon(s)~\doteq~\left(\sup_{x\in \R}\int_{x}^{x+\epsilon}
u_x^2(s,y)\, dy\right)^{1/2}\,.\eeq
Observe that all functions $\eta_\epsilon$ are uniformly bounded.
Moreover, as $\epsilon\to 0$
we have $\eta_\epsilon(t)\downarrow 0$  pointwise at a.e.~time $t$.
Therefore, by the dominated convergence theorem,
\bel{lux}\lim_{\epsilon\to 0}
\int_{\tau}^{t}\left(\int_{x(s)}^{x(s)+\epsilon}
 u_x^2(s,x)dx\right)^{1/2}ds~\leq~\lim_{\epsilon\to 0}
\int_{\tau}^{t}
\eta_\epsilon(s)ds~=~0\,.\eeq
On the other hand, for every time $s\in [\tau,t]$ by construction we have
$$\psi^{\epsilon}_x(s,y)~=~\epsilon^{-1},
\qquad \psi^{\epsilon}_t(s,y)+u(s, x(s))\psi^{\epsilon}_x(s,y)~=~0
\qquad\quad\hbox{for }~ x(s)<y<x(s)+\epsilon\,.$$
This implies
\bel{luxx}\bega{l}\ds\int_{x(s)}^{x(s)+\epsilon}
|\psi^{\epsilon}_t(s,y)+u(s,y)\psi^{\epsilon}_x(s,y)|^2dy
~=~\ds\epsilon^{-2}\int_{x(s)}^{x(s)+\epsilon}
|u(s,y)- u(s, x(s))|^2dy\cr\cr
\leq~\ds \epsilon^{-1}\cdot
\left(\max_{x(s)\leq y\leq x(s)+\epsilon}|u(s,y)- u(s, x(s))|\right)^2
~\leq~\epsilon^{-1}\cdot
\left(\int_{x(s)}^{x(s)+\epsilon} |u_x(s,y)|\, dy\right)^2\cr\cr
\leq~\epsilon^{-1}\cdot
\left(\epsilon^{1/2}\cdot \|u_x(s)\|_{\L^2}\right)^2~\leq~\|u(s)\|_{H^1}\,.
\enda\eeq
Together, (\ref{lux}) and (\ref{luxx}) prove that the integral in (\ref{lexx})
approaches zero as $\epsilon\to 0$.
We now estimate the integral near the corners of the domain:
\bel{corner}
\bega{l}
\ds\left|\left(\int^\tau_{\tau-\epsilon}+\int_t^{t+\epsilon}\right)
\int_{x(s)}
^{x(s)+\kappa\epsilon} u_x(\psi^\epsilon_t+u\psi^\epsilon_x) dx ds\right|\cr\cr
\ds\qquad \leq~\left(\int_{\tau-\epsilon}^{\tau}+\int_t^{t+\epsilon}\right)\left(\int_{x(s)}^{x(s)+\epsilon}
 |u_x|^2dx\right)^{1/2}\left(\int_{x(s)}^{x(s)+\epsilon}
(\psi^\epsilon_t+u\psi^\epsilon_x)^2dx
\right)^{1/2}ds\cr\cr
\ds\qquad \leq~2\epsilon\cdot \|u(s)\|_{H^1}\cdot
\left(\int_{x(s)}^{x(s)+\epsilon} 4\epsilon^{-2}
\left(1+\|u\|_{\L^\infty}\right)^2\, dx
\right)^{1/2}~\leq~C\,\epsilon^{1/2}~\to ~0\enda
\eeq
as $\epsilon\to 0$.
The above analysis has shown that
\[
\lim\limits_{\epsilon\to 0}\int_{\tau-\epsilon}^{t+\epsilon}\int_{x(s)}
^{x(s)+\epsilon} u_x(\psi^\epsilon_t+u\psi^\epsilon_x) dx ds~=~0.
\]
Therefore from (\ref{Inteq}) we recover (\ref{ucar}).
\v
{\bf 6.} Finally, we prove uniqueness of $x(t)$.
Assume there are different $x_1(t)$ and $x_2(t)$, both
satisfying (\ref{ceq})  together with the characteristic equation (\ref{car}).
Choose measurable functions
$\beta_1$ and $\beta_2$ so that
 $x_1(t) = x(t,\beta_1(t))$ and $x_2(t) = x(t,\beta_2(t))$.
Then $\beta_1(\cdot)$ and $\beta_2(\cdot)$ satisfy (\ref{bt}) with the same initial condition. This contradicts with the uniqueness of $\beta$ proved in step {\bf 3}.
\endproof

\v
Relying on (\ref{ucar}) we can now show the Lipschitz continuity of $u$ w.r.t.~$t$, in the
auxiliary coordinate system.
\v
{\bf Lemma 3.} {\it Let $u=u(t,x)$ be a conservative solution of (\ref{1.1}).
Then the  map $(t,\beta)\mapsto u(t,\beta)\doteq u(t, x(t, \beta))$
is Lipschitz continuous, with a constant depending only on the norm $\|u_0\|_{H^1}$.  }
\v
{\bf Proof.}
Using  (\ref{ulip}), (\ref{bt}), and  (\ref{ucar}), and we obtain
$$\bega{rl}
\Big| u(t, x(t,\bar\beta))- u(\tau,\bar\beta)\Big|&\leq~
\Big| u(t, x(t,\bar\beta))- u(t,x(t,\beta(t)))\Big| + \Big| u(t, x(t,\beta(t)))- u(\tau,
x(\tau,\beta(\tau)))\Big| \cr\cr
&\ds\leq~ {1\over2}|\beta(t)-\bar\beta| + (t-\tau)\|P_x\|_{\L^\infty}~\leq~(t-\tau)\cdot (\frac12\|G\|_{\L^\infty}+
\|P_x\|_{\L^\infty}).
\enda
$$\endproof
\v
The next result shows that the solutions $\beta(\cdot)$ of (\ref{bt}) depend Lipschitz
continuously on the initial data.

{\bf Lemma 4.}  {\it Let $u$ be a conservative solution to
(\ref{1.1}). Call $t\mapsto \beta(t;\tau,\bar \beta)$ the solution to the
integral equation
\bel{bt2}
\beta(t)~=~\bar \beta+\int_\tau^t G(\tau, \beta(\tau))\, d\tau,
\eeq
with $G$ as in (\ref{Gdef}).    Then there exists a constant $C$
such that, for any two initial data
$\bar \beta_1,\bar\beta_2$ and any $t,\tau\geq 0$ the corresponding solutions satisfy
\bel{bil}
|\beta(t;\tau,\bar\beta_1)-
\beta(t;\tau,\bar\beta_2)|~\leq~e^{C|t-\tau|}\,|\bar\beta_1-\bar\beta_2|.\eeq
}
\v
{\bf Proof.} Assume $\tau<t$.  By
(\ref{bt2}) it follows
\beq\bega{rl}
|\beta(t;\tau,\bar\beta_1)-\beta(t;\tau,\bar\beta_2)|&=~
\ds\bigg|\bar\beta_1-\bar\beta_2+\int_\tau^t G(s,\beta(s;\tau,\bar\beta_1))-
G(s,\beta(s;\tau,\bar\beta_2))ds\bigg|\cr\cr
&\ds\leq ~|\bar\beta_1-\bar\beta_2|+\int_\tau^t \Big|G(s,\beta(s;\tau,\bar\beta_1))
-G(s,\beta(s;\tau,\bar\beta_2))\Big|\,ds\cr\cr
&\ds\leq ~|\bar\beta_1-\bar\beta_2|+C\int_\tau^t
\Big|\beta(s;\tau,\bar\beta_1)-\beta(s;\tau,\bar\beta_2)\Big|\,ds\cr\cr
&\ds\leq ~|\bar\beta_1-\bar\beta_2|e^{C(t-\tau)},
\enda\eeq
where the last inequality is obtained using
Gronwall's lemma.
The case $t<\tau$ is entirely similar.
\endproof

\v
{\bf Lemma 5.}
{\it Assume $u\in H^1(\R)$ and define the convolution $P$ as in (\ref{Pdef}).
Then $P_x$ is absolutely continuous  and satisfies
\bel{Pxx}
P_{xx}~=~P-\left(u^2+\frac12 u_x^2\right).\eeq
}

{\bf Proof.} The function $\phi(x)=e^{-|x|}/2$ satisfies the distributional identity
$$D^2_x\phi~=~\phi -\delta_0\,,$$
where $\delta_0$ denotes a unit Dirac mass at the origin.
Therefore, for every function $f\in \L^1(\R)$, the convolution satisfies
$$D^2_x(\phi*f)~=~\phi*f - f\,.$$
Choosing $f=u^2 + u_x^2/2$ we obtain the result.
\endproof

\v

\section{Proof of Theorem 2}
\setcounter{equation}{0}
The proof will be worked out in several steps.
\v
{\bf 1.}
By Lemmas 1 and 3,  the map $(t,\beta)\mapsto (x,u)(t,\beta)$ is
Lipschitz continuous.   An entirely similar argument shows that the maps
$\beta\mapsto G(t,\beta)\doteq
G(t,x(t,\beta))$
and $\beta\mapsto P_x(t,\beta)\doteq P_x(t, x(t,\beta))$ are also Lipschitz continuous.
By Rademacher's theorem \cite{Evans}, the partial derivatives
$x_t$,  $x_\beta$, $u_t$, $u_\beta$, $G_\beta$, and $ P_{x,\beta}$ exist almost everywhere.
Moreover, a.e.~point
$(t,\beta)$ is a Lebesgue point for these derivatives.
Calling $t\mapsto \beta(t, \bar \beta)$ the unique solution to
the integral equation (\ref{bt}), by Lemma~4 for a.e.~$\bar \beta$  the following holds.
\begi
\item[{\bf (GC)}] For a.e.~$t> 0$, the point
$(t, \beta(t,\bar\beta))$ is a Lebesgue point for the partial derivatives
$x_t,  x_\beta, u_t, u_\beta, G_\beta, P_{x,\beta}$.  Moreover,
$x_\beta(t, \beta(t,\bar\beta))>0$ for a.e.~$t> 0$.
\endi
If (GC) holds, we then say that $t\mapsto
\beta(t,\bar\beta)$ is a {\bf good characteristic}.
\v
{\bf 2.}
We seek an ODE describing how the quantities
$u_\beta$ and $x_\beta$
vary along a good characteristic.
As in Lemma~4, we denote by $t\mapsto \beta(t;\tau,\bar\beta)$ the solution to
(\ref{bt2}).  If $\tau,t\notin\N$, assuming that $\beta(\cdot; \tau,\bar \beta)$
is a good characteristic, differentiating (\ref{bt2}) w.r.t.~$\bar \beta$ we find
\bel{bode}
{\partial\over\partial\bar\beta} \beta(t;\tau,\bar\beta)~=~
1+\int_\tau^t G_\beta(s,\beta(s;\tau,\bar\beta))
\cdot{\partial\over\partial\bar\beta} \beta(s;\tau,\bar\beta)\,ds
\eeq
Next, differentiating w.r.t.~$\bar\beta$ the identity
$$x(t, \beta(t;\tau,\bar\beta))~=~x(\tau, \bar\beta)+\int_\tau^t
u(s,\,x(s,\beta(t;\tau,\bar\beta)))\, ds$$
 we obtain
\bel{xbode}
x_\beta(t, \beta(t;\tau,\bar\beta))\cdot
{\partial\over\partial\bar\beta} \beta(t;\tau,\bar\beta)~=~x_\beta(\tau,\bar\beta)
+\int_\tau^t u_\beta(s,\beta(s;\tau,\bar\beta))\cdot
{\partial\over\partial\bar\beta} \beta(s;\tau,\bar\beta)\,ds.
\eeq
Finally, differentiating w.r.t.~$\bar\beta$ the identity (\ref{ucar}), we obtain
\bel{ubode}
u_\beta(t, \beta(t;\tau,\bar\beta))\cdot
{\partial\over\partial\bar\beta} \beta(t;\tau,\bar\beta)~=~u_\beta(\tau,\bar\beta)
+\int_\tau^t P_{x,\beta}(s,\beta(s;\tau,\bar\beta))\cdot
{\partial\over\partial\bar\beta} \beta(s;\tau,\bar\beta)\,ds\,.\eeq

Combining (\ref{bode})--(\ref{ubode}),  we thus obtain the system of ODEs
\bel{ODES}\left\{\bega{cl}\ds
{d\over dt}\left[{\partial\over\partial\bar\beta} \beta(t;\tau,\bar\beta)
\right]&\ds=~G_\beta(t,\beta(t;\tau,\bar\beta))
\cdot{\partial\over\partial\bar\beta} \beta(t;\tau,\bar\beta),
\cr
\cr
\ds{d\over dt}\left[x_\beta(t, \beta(t;\tau,\bar\beta))\cdot
{\partial\over\partial\bar\beta} \beta(t;\tau,\bar\beta)
\right]&\ds=~u_\beta(t,\beta(t;\tau,\bar\beta))\cdot
{\partial\over\partial\bar\beta} \beta(t;\tau,\bar\beta),
\cr
\cr
\ds{d\over dt}\left[u_\beta(t, \beta(t;\tau,\bar\beta))\cdot
{\partial\over\partial\bar\beta} \beta(t;\tau,\bar\beta)
\right]&\ds=~P_{x,\beta}(t,\beta(t;\tau,\bar\beta))\cdot
{\partial\over\partial\bar\beta} \beta(t;\tau,\bar\beta).\enda\right.\eeq
In particular, the quantities within square brackets on the left hand sides
of (\ref{ODES}) are absolutely continuous.
From (\ref{ODES}), using Lemma 5 along a good characteristic we obtain
\bel{xubt}\left\{\bega{rl}\ds
{d\over dt}
x_\beta+G_\beta x_\beta& =~u_\beta\,,\cr\cr
\ds{d\over dt}u_\beta+G_\beta u_\beta&\ds =~\left[u^2+\frac12 u_x^2-P\right]x_\beta
~=~\ds\left[u^2+\frac12\left(\frac1{x_\beta}-1\right)-P\right]
x_\beta\cr\cr
&\ds=~\left[u^2-P-\frac12\right]x_\beta+\frac12\,.
\enda\right.\eeq

\v
{\bf 3.} We now go back to the original
$(t,x)$ coordinates and derive an evolution equation
for the partial derivative
$u_x$ along a ``good" characteristic curve.

Fix a point $(\tau,\bar x)$ with $\tau\not\in\N$.
Assume that $\bar x$ is a Lebesgue point for the map  $x\mapsto u_x(\tau, x)$.
Let  $\bar \beta$ be such that $\bar x = x(\tau,\bar\beta)$
and assume that $t\mapsto \beta(t;\tau,\bar\beta)$ is a {\em good characteristic},
so that (GC) holds. We observe that
$$u_x^2(\tau,x)~=~{1\over x_\beta(\tau,\bar\beta)}-1~\geq~0\qquad\qquad x_\beta(\tau,\bar\beta)~>~0\,.$$
As long as $x_\beta>0$, along the characteristic through $(\tau,\bar x)$ the
partial derivative $u_x$ can be computed as
\bel{uxx}
u_x\Big(t,x(t,\beta(t;\tau,\bar\beta))\Big)~=~{u_\beta(t,\beta(t;\tau,\bar\beta))\over
x_\beta(t,\beta(t;\tau,\bar\beta))}\,.\eeq
Using the two ODEs (\ref{xbode})-(\ref{ubode})
describing the evolution of $u_\beta$ and $x_\beta$, we conclude that the map
$t\mapsto u_x(t, x(t,\beta(t;\tau\bar\beta)))$ is absolutely continuous
(as long as $x_\beta\not=0$) and satisfies
\bel{u_xt}\bega{l}\ds
{d\over dt} u_x(t, x(t,\beta(t;\tau\bar\beta)))~=\ds~{d\over dt}\left(
{u_\beta}\over{x_\beta}\right)~=~
\frac{x_\beta\left\{[u^2-P-\frac12]x_\beta+\frac12-u_\beta G_\beta\right\}-u_\beta
\left\{u_\beta-x_\beta G_\beta\right\}}{x_\beta^2}\cr\cr
\qquad \ds=~u^2-P-\frac12+\frac1{2x_\beta}-\frac{u_\beta G_\beta}{x_\beta}-\frac{u_\beta^2}{x_\beta^2}+\frac{u_\beta G_\beta}{x_\beta}~
=~u^2-P-\frac12+\frac1{2x_\beta}-\frac{u_\beta^2}{x_\beta^2}.
\enda\eeq
In turn, as long as $x_\beta>0$ this implies
\bel{arctan}\bega{l}\ds
{d\over dt}\arctan u_x(t, x(t,\beta(t;\tau\bar\beta)))~=~{1\over 1+u_x^2}\cdot {d\over dt} u_x~=~ \left(u^2-P-\frac12+
\frac1{2x_\beta}-\frac{u_\beta^2}{x_\beta^2}\right)x_\beta\cr\cr
\ds\qquad =~ \left(u^2-P-\frac12\right)x_\beta+\frac12-\frac{u_\beta^2}{x_\beta}~=~
\left(u^2-u_x^2-P-\frac12\right)x_\beta+\frac12\cr\cr
\ds\qquad =~ \left(u^2-\frac1{x_\beta}-P+\frac12\right)x_\beta+\frac12~
=~ \left(u^2-P+\frac12\right)x_\beta-\frac12\,.
\enda\eeq
\v
{\bf 4.}
Consider the function
\beq
v~\doteq~\left\{
\bega{cl}
2\arctan u_x \qquad  &\hbox{if}\quad 0<x_\beta\leq 1,\\
\pi \qquad  &\hbox{if}\quad x_\beta=0.
\enda\right.
\eeq
Observe that this implies
\bel{a22}x_\beta~=~{1\over 1+u_x^2}~=~\cos^2{v\over 2}\,,\qquad\quad
{u_x\over 1+u_x^2}~=~{1\over 2}\sin v\,,\qquad\quad {u_x^2
\over 1+u_x^2}~=~\sin^2{v\over 2}\,.\eeq
In the following,  $v$ will be regarded
as a map taking values in the unit circle $\S\doteq [-\pi,\pi]$
with endpoints identified.
We claim that, along each good characteristic, the map $t\mapsto v(t)\doteq
v(t, x(t,\beta(t;\tau\bar\beta)))$ is absolutely continuous and satisfies
\bel{Vt}
{d\over dt}v(t)~=~\left(2u^2-2P+1\right)\cos^2{v\over 2}-1\,.\eeq
Indeed, denote by $x_\beta(t)$, $u_\beta(t)$ and $u_x(t)= u_\beta(t)/x_\beta(t)$
the
values of $x_\beta$, $u_\beta$, and $u_x$ along this particular characteristic.
By (GC) we have $x_\beta(t)>0$ for a.e.~$t>0$.

If $\tau$ is any time where $x_\beta(\tau)>0$,
we can find a neighborhood
$I=[\tau-\delta, \tau+\delta]$ such that  $x_\beta(t)>0$ on $I$.
By (\ref{arctan}) and (\ref{a22}), $v= 2\arctan(u_\beta/x_\beta)$
is absolutely continuous  restricted to $I$ and satisfies (\ref{Vt}).
To prove our claim, it thus remains to show that $t\mapsto v(t)$
is continuous on the null set $\N$ of times where $x_\beta(t)=0$.

Suppose $x_\beta(t_0)=0$. From the identity
\beq
u_x^2(t)~=~{{1-x_\beta(t)}\over{x_\beta(t)}}\,,\eeq
valid as long as $x_\beta>0$, it is clear that $u_x^2\to \infty$ as $t\to t_0$
and $x_\beta(t)\to 0$.   This implies $v(t)=2\arctan u_x(t)
\to \pm\pi$.  Since in $\S$ we identify the points $\pm\pi$, this establishes the
continuity of $v$ for all $t\geq 0$, proving our claim.
\v
{\bf 5.}  Let now $u=u(t,x)$ be a conservative solution.
As shown by the previous analysis, in terms of the variables $t,\beta$ the quantities
$x,u,v$ satisfy the semilinear system

\bel{xuv}\left\{\bega{cl}\ds
{d\over dt} \beta(t, \bar\beta)&=~G(t, \beta(t, \bar\beta)),\cr\cr
\ds {d\over dt} x(t,\beta(t,\bar\beta))&=~u(t,\beta(t,\bar\beta)),\cr\cr
\ds {d\over dt} u(t,\beta(t,\bar\beta))&=~-P_x(t,\beta(t,\bar\beta)),\cr\cr
\ds {d\over dt} v(t,\beta(t,\bar\beta))&=~
\left(2u^2-2P+1\right)\cos^2{v\over 2}-1\,.
\enda\right.\eeq
We recall that $P$ and $G$ were defined at (\ref{Pdef}) and (\ref{Gdef}), respectively.
The function $P$ admits a representation in terms of the variable $\beta$, namely
{\small
\bel{Pc}P(x(\beta))
~= \ds~{1\over 2}
\int_{-\infty}^\infty
\exp\left\{ -\left|
\int_\beta^{\beta'} \cos^2{v(s)\over 2}
\,ds\right|\right\} \cdot
\left[ u^2(\beta')\cos^2{v(\beta')\over 2}+
{1\over 2}\sin^2{v(\beta')\over 2}\right]\,d\beta',
\eeq
\bel{Pxc}P_x(x(\beta))
~= ~{1\over 2}
\left(
\int_\xi^\infty - \int_{-\infty}^\xi\right)
\exp\left\{ -\left|
\int_\xi^{\beta'} \cos^2{v(s)\over 2}\,ds\right|\right\}
\left[ u^2(\beta')\cos^2{v(\beta')\over 2}+
{1\over 2}\sin^2{v(\beta')\over 2}\right]\,d\beta'.
\eeq
}
For every $\bar\beta\in \R$ we have the initial condition
\bel{ico}\left\{\bega{cl}\ds\beta(0,\bar\beta)&=~\bar\beta,\cr
x(0,\bar\beta)&=~x(0,\bar\beta),\cr
u(0,\bar\beta)&=~u_0(x(0,\bar\beta)),\cr
v(0,\bar\beta)&=~2\arctan u_{0,x}(x(0,\bar\beta)).\enda\right.
\eeq
By the Lipschitz continuity of all coefficients,
the Cauchy problem (\ref{xuv}), (\ref{ico}) has a unique
solution, globally defined for all $t\geq 0$, $x\in \R$.
\v
{\bf 6.}
To complete the proof of uniqueness,   consider two conservative solutions
$u,\tilde u$ of the Camassa-Holm equation (\ref{1.1})
with the same initial data $u_0\in H^1(\R)$.
For a.e.~$t\geq 0$  the
corresponding  Lipschitz continuous
maps $\beta\mapsto x(t,\beta)$, $\beta\mapsto \tilde x(t,\beta)$
are strictly increasing.  Hence they have continuous inverses, say
$x\mapsto \beta^* (t,x)$, $x\mapsto \tilde \beta^* (t,x)$.

By the previous analysis, the map    $(t,\beta)\mapsto (x,u,v)(t,\beta)$
is uniquely determined by the initial data $u_0$.  Therefore
$$x(t,\beta)~=~\tilde x(t,\beta),\qquad\qquad u(t,\beta)=\tilde u(t,\beta).$$
In turn, for a.e.~$t\geq 0$ this implies
$$u(t,x) ~=~u(t,\beta^*(t,x))~=~\tilde u(t,\tilde \beta^*(t,x))~=~\tilde u(t,x).$$
\endproof

\end{document}